\documentclass [11pt]{article}
\usepackage{amsmath, amsthm, amssymb}
\usepackage{graphicx}
\newfont{\bbb} {msbm10}
\newcommand{\R}{\Bbb{R}}
\newcommand{\bS}{\Bbb{S}}

\newcommand{\sbs}{\subset}
\newcommand{\ra}{\rightarrow}

\newcommand{\hg}{{\hat{g}}}
\newcommand{\hh}{{\hat{h}}}
\newcommand{\hf}{{\hat{f}}}

\newcommand{\p}{\partial}

\newcommand{\cE}{{\cal{E}}}

\newcommand{\HH}{\Bbb{H}}

\newcommand{\ssl}{{_\lambda}}

\newcommand{\0}[1]{_{_{#1}}}
\newcommand{\rC}{{\rm{C}}\,}

\usepackage[all]{xy}

\begin{document}

\title{Cut Limits on Hyperbolic Extensions}
\author{Pedro Ontaneda\thanks{The author was
partially supported by a NSF grant.}}
\date{}

\maketitle

\begin{abstract} Hyperbolic extensions were defined and studied in \cite{O3}.
Cut limits of families of metrics were introduced in \cite{O4}.
In this paper we show that if a family of metrics $\{h\0{\lambda}\}$ has cut limits then
the family of hyperbolic extensions $\{ \cE_k(h\0{\lambda})\}$ also has cut limits.

The results in this paper are used in the problem of smoothing
Charney-Davis strict hyperbolizations \cite{ChD}, \cite{O}.
\end{abstract}
\vspace{.3in}

\noindent {\bf \large  Section 1. Introduction.}

This paper deals with the relationship between two
concepts:  ``hyperbolic extensions", which were studied in \cite{O3}, and
``cut limits of families of metrics", which were defined in \cite{O4}. Before stating our main result we first introduce
these concepts here.\\

\noindent {\bf A. Hyperbolic extensions.}
 Recall that the hyperbolic $n$-space $\HH^n$ is isometric to $\HH^{k}\times \HH^{n-k}$ with
warp product metric $(\cosh ^2\, r)\,\sigma\0{\HH^{k}}+\sigma\0{\HH^{n-k}}$, where $\sigma\0{\HH^{l}}$ denotes the hyperbolic metric of
$\HH^{l}$, and  $r:\HH^{n-k}\ra[0,\infty)$ is the distance to a fixed point in $\HH^{n-k}$.
For instance, in the case $n=2$, since $\HH^1=\R^1$ we have that $\HH^2$ is isometric to $\R^2=\{(u,v)\}$ with metric $\cosh ^2v\, du^2+ dv^2$. 
The concept of ``hyperbolic extension" is a generalization of
this construction; we explain this in the next paragraph.\\

Let $(M^n,h)$ be a complete Riemannian manifold with {\it center} $o=o_{_M}\in M$, that is, the exponential map $exp_o:T_oM\ra M$ is a diffeomorphism.
The warp product metric \vspace{.1in}

\hspace{2in}$f=(\cosh ^2 r)\, \sigma\0{\HH^k}+h.$\vspace{.1in}

\noindent on $\HH^k\times M$
is the {\it hyperbolic extension (of dimension $k$)} 
of the metric $h$. Here $r$
is the distance-to-$o$ function on $M$.
We write $\cE_k(M)=(\HH^k\times M,f)$, and $f=\cE_k(h)$.
We also say that $\cE_k(M)$ is the {\it hyperbolic extension
(of dimension $k$) of $(M,h)$} (or just of $M$).
Hence, for instance, we have $\cE_k(\HH^l)=\HH^{k+l}$.  Also write $\HH^k=\HH^k\times \{o\0{M}\}\sbs\cE_k(M)$ and
we have that any $p\in\HH^k$ is a center of $\cE_k(M)$
(see Remark 2.3 (3)). \\

\noindent {\bf Remarks 1.1.} \\
{\bf 1.} Let $M^n$ have center $o$.
Using a fixed orthonormal basis on $T_oM$ and the exponential
map we can identify $M$ with $\R^n$, and $M-\{o\}$ with
$\R^n-\{0\}=\bS^{n-1}\times (0,\infty)$. Hence the spheres
$\bS^{n-1}\times \{r\}\sbs\bS^{n-1}\times (0,\infty)$ are
geodesic spheres, and the rays $t\mapsto tv= ( v , t)\in\bS^{n-1}
\times (0,\infty)=M-\{o\}$, are distance minimizing 
geodesic rays emanating from
the center.\\
{\bf 2.} Let $g'$ be another metric on $M$. Suppose we
can write $g'=g'_r+dr^2$ on $\bS^{n-1}\times (0,\infty)=
M-\{o\}$, (this last identification is done using $g$).
Then the geodesic spheres around $o$, and the
geodesics emanating from $o$ for $g$ and $g'$
coincide.\\

\noindent {\bf B.  Cut limits.}
Before we talk about ``cut limits" we need some preliminary
definitions and facts. 
Let $(M^n,g)$ have center $o$.
Then the metric $g$ (outside the center) has the form $g=g_r+dr^2$. Here we are identifying (see 1.1) the space
$M-\{o\}$ with $\bS^{n-1}\times (0,\infty)$, thus each
$g_r$ is a metric on the sphere $\bS^{n-1}$.\\

\noindent {\bf Examples.}\\
{\bf 1.} The Euclidean metric $\sigma\0{\R^n}$ on
$\R^n$ can be written as $\sigma\0{\R^n}=g_r+dr^2$
with $g_r=r^2 \sigma\0{\bS^{n-1}}$, where
$\sigma\0{\bS^{n-1}}$ is the round metric
on the sphere $\bS^{n-1}$.\\
{\bf 2.} The hyperbolic metric $\sigma\0{\HH^n}$ on
$\R^n$ can be written as $\sigma\0{\HH^n}=g_r+dr^2$
with $g_r=\sinh^2(r)\, \sigma\0{\bS^{n-1}}$.\\

Let $(M,g)$ have center $o$ and write  $g=g_r+dr^2$.
Let $r\0{0}>0$. We can think of the metric $g\0{r\0{0}}$
as being obtained from $g=g_r+dr^2$ by ``cutting" $g$ along
the sphere of radius $r\0{0}$, so we call
the metric $g\0{r\0{0}}$ on
$\bS^{n-1}$ the {\it spherical cut of
$g$ at $r\0{0}$.} Let\vspace{.1in}

\noindent {\bf (1.2)}$\hspace{2.2in}\hat{g}\0{r\0{0}}
=\big(  \frac{1}{\sinh^2(r\0{0})}\big) g\0{r\0{0}}.$\vspace{.13in}

We call the metric $\hat{g}\0{r\0{0}}$ on
$\bS^{n-1}$ given by (1.2) the {\it normalized spherical
cut of $g$ at $r\0{0}$}. In the particular case that $g=g_r+dr^2$ is a warped-by-$\sinh$
metric we have $g_r=\sinh^2(r)g'$ for some fixed
$g'$ independent of $r$. In this case the spherical cut
of $g=\sinh^2(r)g'+dt^2$ at $r\0{0}$ is $\sinh^2(r\0{0})g'$,
and the the normalized spherical cut at $r\0{0}$ is
$\hat{g}\0{r\0{0}}=g'$.\\

\noindent {\bf Example.} If $g=\sigma\0{\HH^n}=\sinh^2(r)\, \sigma\0{\bS^{n-1}}+dr^2$,
the normalized spherical cut at $r\0{0}$ is $({\widehat{\sigma\0{\HH^n}}})\0{r\0{0}}=\sigma\0{\bS^{n-1}}$. And the  spherical
cut at $r\0{0}$ is $\sinh^2(r\0{0})\, \sigma\0{\bS^{n-1}}$.\\

Let $(M^n,g)$ have center $o$.
We now consider families of metrics $\{ g_\ssl\}_{\lambda> \lambda_0}$ on $M$ of the form $g_\ssl=\big(g_\ssl\big)_r+dr^2$. Here $\lambda_0 >0$, and the identification
$M-\{o\}=\bS^{n-1}\times (0,\infty)$ is done using $g$; see Remark 1.1. 
We call such a family an {\it $\odot$-family of metrics on $(M,g)$}. 
(We use the symbol $\odot$ to evoke the idea that all metrics $g\0{\lambda}$ have a common center and spheres).
The reason we are interested in  these
families is that they are key ingredients
in Riemannian Hyperbolization \cite{O} (also see \cite{O4}).
Moreover, the Main Theorem in this paper is used in \cite{O}.\\

Let $b\in \R$. By cutting each $g_\ssl$ at 
$b+\lambda$ we obtain a one-parameter family
$\{\widehat{\big(g_\ssl\big)}\0{\lambda+b}\}\0{\lambda}$
of metrics on the sphere $\bS^{n-1}$.
(The metric ${\widehat{\big(g_\ssl\big)}}_{_{\lambda+b}}$ is the normalized spherical cut of $g_\ssl$ at
$\lambda+b$). Here $\lambda>max\{\lambda_0,-b\}$, % if $-b\geq\lambda_0$, 
so that the definition makes sense.
We say that the $\{g_\ssl\}$ has {\it cut limit at $b$} if this family $C^2$-converges, as $\lambda\ra\infty$. That is, there is a 
$C^2$ metric $\hg_{_{\infty}}^b$ on $\bS^{n-1}$ such that \vspace{.1in}

\noindent {\bf (1.3)}\hspace{1.5in}$
\big|\,\,{\widehat{\big(g_\ssl\big)}}_{_{\lambda+b}} \,\,-\,\, \hg_{_{\infty}}^b\,\,\big|_{C^2(\bS^{n-1})}\,\,\longrightarrow\,\, 0
$  \,\,\,\,\,\,\,as\,\,\,\,\,\, $\lambda\ra \infty.$\vspace{.1in}

Here the arrow means convergence in the $C^2$-norm on the 
space of $C^2$ metrics on  $\bS^{n-1}$.

\noindent {\bf Remark 1.4.}
The  $C^2$ norm is taken with respect to a 
fixed locally finite atlas with
extendable charts, i.e. charts that can be 
extended to the (compact) closure of their domains.\\

Let $I\sbs\R$ be an interval (compact or noncompact). 
We say that the $\odot$-family $\{g_\ssl\}$ has 
{\it cut limits on $I$} if the convergence in (1.3) is 
 uniform with compact supports in the variable
in $b\in I$. Explicitly this means: for every $\epsilon>0$,
and compact $K\sbs  I$ there is $\lambda_*$  such that
{\scriptsize $\big|\,\,{\widehat{\big(g_\ssl\big)}}_{_{\lambda+b'}} \,\,-\,\, \hg_{_{\infty+b'}}\,\,\big|\0{C^2(\bS^{n-1})}<\epsilon$}, for $\lambda>\lambda_*$ and $b'\in K$.\\

\noindent {\bf Remark 1.5.} Equivalently, the $\odot$-family $\{g_\ssl\}$ has 
{\it cut limits on $I$} if
for every $\epsilon>0$,
and $b\in I$ there are $\lambda_*$ and neighborhood
$U$ of $b$ in $I$ such that
{\scriptsize $\big|\,\,{\widehat{\big(g_\ssl\big)}}_{_{\lambda+b'}} \,\,-\,\, \hg_{_{\infty+b'}}\,\,\big|\0{C^2(\bS^{n-1})}<\epsilon$}, for $\lambda>\lambda_*$ and $b'\in U$.\\

If $\{g_\ssl\}$ has cut limits on $I$ then it 
has a cut limit at $b$, 
for every $b\in I$. Finally we say that 
$\{g_\ssl\}$ has a {\it cut limits} if
$\{g_\ssl\}$ has a cut limits on $\R$.\\

\noindent {\bf Remark 1.6.} If $\{g_\ssl\}\0{\lambda}$ is a family of metrics then
$\{g\0{\lambda(\lambda')}\}\0{\lambda'}$ is a {\it reparametrization} of 
$\{g_\ssl\}\0{\lambda}$, where $\lambda'\mapsto\lambda(\lambda')$
is a change of variables. For instance, if we use translations, the following
holds: $\{g_\ssl\}\0{\lambda}$ has cut limits at $b$ if and only if
$\{g\0{\lambda'+a}\}\0{\lambda'}$ has cut limits at $b+a$; here
the change of variables is $\lambda=\lambda'+a$. \\

\noindent {\bf C. Statement of main result.}
 Here is a natural question:\vspace{.1in}

\noindent {\bf Question.} {\it If
$\{h_\ssl\}\0{\lambda}$ has cut limits, does 
$\{\cE_k(h_\ssl)\}\0{\lambda}$ have cut limits?} \vspace{.1in}

\noindent {\bf Remark.} 
More generally we can ask whether $\{\cE_k(h\0{\lambda})\}\0{\lambda'}$
has cut limits, where $\lambda=\lambda(\lambda')$. Of course the answer
would depend on the change of variables $\lambda=\lambda(\lambda')$.\\
 
Our main result gives an affirmative answer to this
question provided the family $\{h_\ssl\}$ 
is, in some sense, nice near the origin. Explicitly,
we say that $\{h_\ssl\}\0{\lambda>\lambda\0{0}}$ is {\it hyperbolic
around the origin} if there is a $B\in \R$ such that
\vspace{.05in}

\hspace{2.3in}
${\widehat{\big(h_\ssl\big)}}_{_{\lambda+b}}=
\sigma\0{\bS^{n-1}}.$\vspace{.1in}

\noindent for every $b\leq B$ and every 
$\lambda>max\{\lambda_0,-b\}$.
Note that this implies that each $h_\ssl$ is
canonically hyperbolic on the ball of radius $\lambda +B$,
i.e. $h_\ssl=\sinh^2(r)\,\sigma\0{\bS^{n-1}}+dr^2$ on the ball
of radius $\lambda +B$.
Examples of $\odot$-families that are
hyperbolic around the origin are families obtained
using hyperbolic forcing \cite{O4}.\\

As mentioned before our main result answers affirmatively the question above.
Moreover it also says that some reparametrized 
families $\{\cE_k(h\0{\lambda})\}\0{\lambda'}$ 
have cut limits as well,
for certain change of variables $\lambda=\lambda(\lambda')$.
Write $\lambda=\lambda(\lambda',\theta)=
\sinh^{-1}(\sinh(\lambda')\,\sin\,\theta)$, for a fixed $\theta$.
We say that $\{\cE_k(h\0{\lambda})\}\0{\lambda'}$ 
is the {\it $\theta$-reparametrization
of} $\{\cE_k(h\0{\lambda})\}\0{\lambda}$. 
Note that if we
consider an hyperbolic right triangle
with one angle equal to $\theta$ and 
side (opposite to $\theta$) of length
$\lambda$, then 
$\lambda'$ is the
length of the hypotenuse of the triangle. 
All $\theta$-reparametrizations, in the limit $\lambda'\ra \infty$, differ
just by translations; that is, a simple calculation shows
that $\lim_{\lambda'\ra \infty}\lambda(\lambda')
-\lambda'=\ln\sin\theta$.
We are now ready to state our Main result.\\

\noindent {\bf  Main Theorem.} 
{\it Let $M$ have center $o$. Let  $\{h_\ssl\}\0{\lambda}$ be $\odot$-family of metrics
on $M$.
If $\{h_\ssl\}\0{\lambda}$ is 
hyperbolic around the origin and  has cut limits, 
then for every $\theta\in(0,\pi/2]$ the $\theta$-reparametrization 
$\{\cE_k(h\0{\lambda})\}\0{\lambda'}$ has cut limits.}\\

Note that $\theta=\pi/2$ gives $\lambda=\lambda'$
answering the question above. The paper is structured as follows.
In Section 2 we review some facts about
hyperbolic extensions. In Section 3
we introduce useful coordinates on the
spheres of a hyperbolic extension.
In Section 4 we study normalized spherical cuts on
hyperbolic extensions.
Finally in Section 5 we deal with cut limits
in a bit more detail and prove the Main Theorem.\vspace{.3in}

\noindent {\bf \large Section 2.  Hyperbolic Extensions.}

Notational convention: we will denote all fixed centers on manifolds by the same letter ``$o$". If the manifold $M$ needs to be specified
we will write $o=o\0{M}$, which means that $o$ is a center in $M$.\vspace{.1in}

Note that $\HH^k$ is convex in 
$\cE_k(M)$ (see \cite{BisOn}, p.23).
Let $\eta$ be a complete geodesic in $M$ passing though $o$
and let $\eta^+$ be one of its two geodesic rays (beginning at $o$) . Then $\eta$ is 
a totally geodesic subspace of $M$ and $\eta^+$ is 
convex (see \cite{O3}). 
Also, let $\gamma$ be a complete geodesic in $\HH^k$. 
The following two results are proved in Section 3 of \cite{O3}.
\vspace{.1in}

\noindent {\bf Lemma 2.1.} {\it 
The subspace\,\, $\gamma\times \eta^+$ is a convex in $\cE_k(M)$,
and $\gamma\times \eta$ is totally geodesic in $\cE_k(M)$.}
\vspace{.1in}

\noindent {\bf Corollary 2.2.} {\it The subspaces 
\,$\HH^k\times\eta^+$  and $\gamma\times M$ 
are convex in $\cE_k(M)$. Also  \,$\HH^k\times\eta$
is totally geodesic in $\cE_k(M)$.}\vspace{.1in}

\noindent {\bf Remarks 2.3.}

\noindent {\bf 1.} By convexity above we mean the following:
a set $A$ is convex if given two points in $A$ 
any distance minimizing geodesic joining these points lies in
$A$.

\noindent {\bf 2.} As pointed out in Section 3 of \cite{O3},
the proof of Lemma 2.1 (which is Lemma 3.1 in \cite{O3})
can easily be adapted to show that $\{y\}\times M$ are convex in 
$\cE_k(M)$. Alternatively, it is not hard to prove that $\{y\}\times M$
is convex in $\gamma\times M$; this together with Corollary 3.2 implies that
$\{y\}\times M$ are convex in 
$\cE_k(M)$.

\noindent {\bf 3.} Note that $\HH^k\times\eta$ (with metric
induced by $\cE_k(M)$)
is isometric to $\HH^k\times \R$ with warp product metric 
$\cosh ^2 v\, \sigma\0{\HH^k}+dv^2$, which is just hyperbolic 
$(k+1)$-space $\HH^{k+1}$. Also $\gamma\times\eta$
is isometric to $\R\times \R$ with warp product metric 
$\cosh ^2 v\, du^2+dv^2$, which is just hyperbolic 2-space 
$\HH^2$.
In particular every point in 
$\HH^k=\HH^k\times\{ o\} \sbs\cE_k(M)$ is a center point.\\

As before we use $h$ to identify $M-\{ o\}$ with $\bS^{n-1}\times \R^+$. Sometimes we will denote a point
$v=(u,r)\in \bS^{n-1}\times\R^+=M-\{ o\}$ by $v=ru$.
Fix a center $o\in \HH^k\in \cE_k(M)$.  Then, for $y\in\HH^k-\{ o\}$ we can also write $y=t\,w$, $(w,t)\in \bS^{k-1}\times\R^+$. 
Similarly, using the exponential map we can identify $\cE_k(M)-\{ o\}$
with $\bS^{k+n-1}\times \R^+$, and for $p\in\cE_k(M)-\{ o\}$
we can write $p=s\,x$, $(x,s)\in\bS^{k+n-1}\times\R^+$.\vspace{.1in}

We denote the metric on $\cE_k(M)$ by $f$ and we can write $f=f_s+ds^2$.
Since $\HH^k$ is convex in $\cE_k(M)$ we can write $\HH^k-\{ o\}=\bS^{k-1}\times \R^+\sbs\bS^{k+n-1}\times \R^+$
and $\bS^{k-1}\sbs \bS^{k+n-1}$.\vspace{.1in}

A point $p\in\cE_k(M)\, -\, \HH^k$ has two sets of coordinates: the {\it polar coordinates}
$(x,s)=(x(p),s(p))\in \bS^{k+n-1}\times \R^+$ and the {\it hyperbolic extension coordinates} $(y,v)=(y(p), v(p))\in \HH^k\times M$. Write $M_o=\{o\}\times M$.
Therefore we have the following functions:
{\small $$
\begin{array}{lll}
{\mbox{the distance to {\it o} function:}}  
& s:\cE_k(M)\ra [0,\infty), & s(p)=d\0{\cE_k(M)}(p,o)\\\\
{\mbox{the direction of {\it p} function:}}  
& x:\cE_k(M)-\{o\}\ra \bS^{n+k-1}, & p=s(p)x(p)\\\\
{\mbox{the distance to {\it $\HH^k$} function:}}  
& r:\cE_k(M)\ra [0,\infty), & r(p)=d\0{\cE_k(M)}(p,\HH^k)\\\\
{\mbox{the projection on $\HH^k$ function:}}  
& y:\cE_k(M)\ra \HH^k, &\\\\
{\mbox{the projection on $M$ function:}}  
& v:\cE_k(M)\ra M, & \\\\
{\mbox{the projection on $\bS^{n-1}$ function:}}  
& u:\cE_k(M)-\HH^k\ra \bS^{n-1}, & v(p)=r(p)u(p)\\\\
{\mbox{the length of $y$ function:}}  
& t:\cE_k(M)\ra [0,\infty), & t(p)=d_{\HH^k}(y(p),o)\\\\
{\mbox{the direction of $y$ function:}}  
& w:\cE_k(M)-M_o\ra \bS^{k-1}, & y(p)=t(p) w(p)
\end{array}
$$}

Note that $r=d_M(v, o)$. Note also that, by Lemma 2.1, the functions $w$ and $u$ are constant on geodesics emanating from $o\in\cE_k(M)$, that is
$w(sx)=w(x)$ and $u(sx)=u(x)$.\vspace{.1in}

Let $\p_r$ and $\p_s$ be the gradient vector fields of $r$ and $s$, respectively. Since the $M$-fibers $M_y=\{ y\}\times M$ are convex
the vectors $\p_r$ are the velocity vectors of the speed one geodesics of the form $a\mapsto (y, a\, u)$, $u\in\bS^{n-1}\sbs M$. These geodesics 
emanate from (and  orthogonally to) $\HH^k\sbs \cE_k(M)$.
Also the vectors  $\p_s$ are the velocity vectors of the speed one geodesics 
emanating from $o\in\cE_k(M)$. For $p\in\cE_k(M)$, denote by $\bigtriangleup =\bigtriangleup (p)$ the right triangle with vertices $o$, $y=y(p)$, $p$
and sides the geodesic segments $[o,p]\in\cE_k(M)$, $[o,y]\in\HH^k$, $[p,y]\in\{ y\}\times M\sbs\cE_k(M)$.
(These geodesic segments are unique and well-defined because:\, (1) $\HH^k$ is
convex in $\cE_k(M)$,\, (2) $(y,o)=o_{_{\{ y\}\times M}}$ and $o$ are centers in $\{ y\}\times M$ and $\HH^k\sbs\cE_k(M)$, respectively.)\vspace{.1in}

Let $\alpha:\cE_k(M)-\HH^k\ra \R$ be the angle between 
$\p_s$ and $\p_r$ (in that order), thus  $\cos\, \alpha=f(\p_r,\p_s)$, $\alpha\in [0,\pi]$. 
Then $\alpha=\alpha(p)$ is the interior angle, at $p=(y,v)$, of the right triangle $\bigtriangleup =\bigtriangleup (p)$.
We call  $\beta(p)$ the interior angle of this triangle at $o$, that is $\beta(p)=\beta(x)$ is the spherical distance 
between $x\in \bS^{k+n-1}$ and the totally geodesic sub-sphere $\bS^{k-1}$. Alternatively, $\beta$ is the angle between the geodesic segment
$[o,p]\sbs\cE_k(M)$ and the convex submanifold $\HH^k$.
Therefore $\beta$ is constant on each geodesic emanating from $o\in\cE_k(M)$, that is
$\beta(sx)=\beta(x)$. The following corollary follows from Lemma 2.1 (see 3.1 in \cite{O3}).\vspace{.1in}

\noindent {\bf Corollary 2.4.} {\it Let $\eta^+$  (or $\eta$) be a geodesic ray (line) in $M$ through $o$ containing
$v=v(p)$ and $\gamma$ a geodesic line in $\HH^k$ through $o$ containing $y=y(p)$. 
Then $\bigtriangleup (p)\sbs \gamma\times \eta^+\sbs 
\gamma\times \eta$.}\\

\noindent {\bf Remark 2.5.}
Note that the right geodesic triangle $\bigtriangleup (p)$ 
has sides of length $r=r(p)$, $t=t(p)$ and $s=s(p)$. 
By Lemma 2.1 and Remark 2.3(3)
we can consider $\bigtriangleup$ as contained in a 
totally geodesic copy of hyperbolic 2-space $\HH^2(p)$.
The plane $\HH^2(p)$ is well defined for $p$
outside $\HH^k\cup (\{o \}\times M)$. We will write
$\HH^2(p)=\gamma_w\times\eta_u$, where $p=(y,v)\in
\HH^k\times M$, $y=tw$, $v=ru$.\vspace{.1in}

Hence, by Remark 2.5, using hyperbolic trigonometric identities
we can find relations between $r$, $t$, $s$, $\alpha$ and $\beta$. For instance, using the hyperbolic law of sines we get:

\vspace{.1in}

\noindent {\bf (2.6)}\hspace{1.9in}
$\sinh\, (r)\, =\, \sin\, (\beta)\,\, \sinh\, (s)$\vspace{.1in}

In Section 4 we will need the following result.
\vspace{.1in}

\noindent {\bf Proposition 2.7.} {\it  The following identity holds outside $\HH^k\cup \big(\{o \}\times M$\big)}

$$
\big( \sinh^2 (s) \big)\, d\beta\,^2\,\,\, +\,\,\, ds^2\,\,\, =\,\,\, \cosh ^2 (r) \, dt^2\,\,\,+\,\,\, dr^2
$$

\noindent {\bf Proof.} First a particular case. Take $M=\R$ and $k=1$, hence $\cE_k(M)=\cE_1(\R)=\HH^2$. In this case the left-hand side of
the identity above is the expression of the metric of $\HH^2$ in polar coordinates $(\beta, s)$, and right hand side of the equation is
the expression of the same metric in the hyperbolic extension coordinates $(r,t)=(v,y)$. (Here $r$ and $t$ are ``signed" distances.)
Hence the equation holds in this particular case.  This particular case, together with the fact that
$\HH^2(p)$ is isometric to $\HH^2$, and the 
following claim prove the proposition.\vspace{.1in}

\noindent {\bf Claim.} {\it The
functionals $d\beta$, $ds$, $dt$, $dr$, at $p\in \HH^k\cup (\{o \}\times M)$, are zero on vectors perpendicular to 
$\HH^2(p)$.}
\vspace{.1in}

\noindent {\bf Proof of the Claim}.
To prove the claim let $u$ be a vector perpendicular to $\HH^2(p)$, at $p$. Since the ray
$s\mapsto sx(p)$ is contained in $\HH^2(p)$ we have that $u$ is tangent to the
sphere of radius $s(p)$ centered at $o$. Therefore $ds(u)=0$. 
\vspace{.1in}

Next we prove that
$dr(u)=0$ and $dt(u)=0$. Note that $u$ is a linear combination of vectors
perpendicular to  $\HH^2(p)$ that are either tangent to $\{y\}\times M$ or
$\HH^k\times\{v\}$, where $y=y(p)$ and $v=v(p)$. Therefore it is enough
to assume $u$ is tangent to $\{y\}\times M$ or
$\HH^k\times\{v\}$.\vspace{.1in}

First assume
that  $u$ is perpendicular to  $\HH^2(p)$ and tangent to $\{y\}\times M$.
Since $u$ is tangent to $\{y\}\times M$ we get that $dt(u)$. And since
$u$ is perpendicular to the ray $r\mapsto rv$ in $\{y\}\times M$
(because this ray is contained in $\HH^2(p)$)
we get that $dr(u)=0$.\vspace{.1in}

Next assume that  $u$ is perpendicular to  $\HH^2(p)$ and tangent to
$\HH^k\times\{v\}$. Then $dr(u)=0$. And since
$u$ is perpendicular to the ray $t\mapsto ty$ in $\HH^k\times\{v\}$
(because this ray is contained in $\HH^2(p)$)
we get that $dt(u)=0$. \vspace{.1in}

Finally, the equation $d\beta(u)=0$ follows from $ds(u)=0$, $dt(u)=0$, $dr(u)=0$,
the fact that $\beta$ is a function of $s,\, t,\, r$, and the chain rule.
This proves the claim and concludes the
proof of Proposition 2.7.\\
\vspace{.2in}

\noindent {\bf \large Section 3. Coordinates On The Spheres $S_s\big( \cE_k(M)\big)$.} 

Let $N^n$ have center $o$. The geodesic sphere of radius $r$ centered at $o$ will be denoted by $\bS_r=\bS_r(N)$ and we can identify $\bS_r$ with $\bS^{n-1}\times\{ r\}$.\vspace{.1in}

Let $M$ have center $o$ and metric $h$.  Consider the hyperbolic extension $\cE_k(M)$ of $M$
with center $o\in\HH^k=\HH^k\times \{ o\}\sbs \cE_k(M)$ and metric $f$.
Since $\HH^k\sbs \cE_k(M)$ is convex, we can write $\bS_s(\cE_k(M))\cap\HH^k=\bS_s(\HH^k)$. Equivalently 
$\big(\bS^{k+n-1}\times\{s\}\big)\cap \HH^k=\bS^{k-1}\times\{s\}$.
Write $M_o=\{ o\}\times M$. Also write 
\begin{center}$
E_k(M)\, =\, \cE_k(M)\, -\, \big( \HH^k\, \coprod \, M_o  \big)
$\end{center}and
\begin{center}$
S_s\big(\, \cE_k(M)\, \big)\, =\,\bS_s\big( \cE_k(M) \big)\, \bigcap\, E_k(M)\,=
\, \bS_s\big(\, \cE_k(M)\, \big)\, -\, \big( \HH^k\, \coprod \, M_o  \big)
$\end{center}
Note that the functions $\alpha$ and $\beta$ are well-defined and smooth on $E_k(M)$, and $0<\beta(p)<\pi/2$.
Moreover, by Remark 2.5, the plane $\HH^2(p)=\gamma_w\times \eta_u$ is well defined for $p\in E_k(M)$. 
As in Remark 2.5, here $p=(y,v)\in
\HH^k\times M$, $y=tw$, $v=ru$. 
Recall that $\bigtriangleup (p)\sbs\HH^2(p)$ (see
Corollary 2.4 and Remark 2.5). \vspace{.1in}

By the identification between $\bS^{n+k-1}\times\{s\}$ with $\bS_s(\cE_k(M))$
and Lemma 2.1 we have that $\HH^2(p)\cap\,\bS_s(\cE_k(M))$ gets identified with a
geodesic circle $\bS^1(p)\sbs\bS^{n+k-1}$. 
Moreover, since $\HH^2(p)$ and $\HH^k$ intersect orthogonally on $\gamma_w$, we have that the spherical geodesic segment 
$\big[ x(p), w(p)\big]\0{\bS^{n+k-1}}$
intersects $\bS^{k-1}\sbs\bS^{n+k-1}$ orthogonally at $w$. This together with the fact that $\beta<\pi/2$ imply that $\big[ x(p), w(p)\big]_{\bS^{n+k-1}}$ is a length minimizing spherical geodesic in $\bS^{k+n-1}$
joining $x$ to $w$. Consequently $\beta=\beta(p)$ is the length of $\big[ x(p), w(p)\big]\0{\bS^{n+k-1}}$.\vspace{.1in}

We now give a set of coordinates on $S_s\big(\cE_k(M)\big)$. For 
$p\in S_s\big(\cE_k(M)\big)$  define 
\begin{center}$
\Xi\, (p)\,=\, \Xi_s(p)\,=\,\big(\,  w\, ,\,  u  \, ,\,  \beta  \,   \big)\,\in\, \bS^{k-1}\times\, \bS^{n-1}\times \big( 0\, ,\, \pi/2  \big)
$\end{center}
\noindent where $w=w(p)$, $u=u(p)$, $\beta=\beta(p)$.
Note that $\Xi$ is constant on geodesics emanating from $o\in\cE_k(M)$, that is $\Xi(sx)=\Xi(x)$.\vspace{.1in}

Using hyperbolic trigonometric identities (e.g. identity 2.6) we can find well defined and smooth functions $r=r(s,\beta)$ and $t=t(s,\beta)$ such that
$r$, $s$, $t$ are the lengths of the sides of a right geodesic triangle on $\HH^2$ with angle $\beta$ opposite the the side with length $r$.
With these functions we can construct explicitly a smooth inverse to $\Xi$.\vspace{.1in}

\noindent {\bf 3.1 Remarks. } \\
{\bf 1.} For $(w,u)\in \bS^{k-1}\times\bS^{n-1}$ we have 
\begin{center}$\Xi\bigg(  \big(\gamma_w\times\eta_u\big) \,\cap\, S_s\big(\cE_k(M)\big)\bigg)\,=\,\{ \pm w\}\times\{\pm u\}\times (0,\pi/2)$\end{center}
\noindent By Lemma 2.1 the paths $a\mapsto (\pm w,\pm u, a)$ four spherical (open) geodesic segments emanating orthogonally from $\bS^{k-1}$. \\
{\bf 2.} For $w\in \bS^{k-1}$ we have 
\begin{center}$\Xi\bigg(  \big(\gamma_w\times M\big)\,\cap\, S_s\big(\cE_k(M)\big) \bigg)\,=\,\{\pm w\}\times\bS^{n-1}\times (0,\pi/2)$\end{center}
\noindent By Corollary 2.2 we have that this set is a spherical geodesic ball of radius $\pi/2$ and of dimension $n$ (with its center deleted)
intersecting  $\bS^{k-1}$ orthogonally at $w$. Note that the geodesic segments on this ball emanating from $w$ are the spherical
geodesic segments of item 1, for all $u\in \bS^{n-1}$.\\
{\bf 3.} For $w\in \bS^{k-1}$ and $r$ with $0<r< s$ we have 
\begin{center}$\Xi\bigg(\big(  \gamma_w\times \bS_r(M)\big) \cap\, S_s\big(\cE_k(M)\bigg)\, =\, \{ w\}\times\bS^{n-1}\times \beta(r)$\end{center}
\noindent where $\beta(r)$ is the angle of the right geodesic hyperbolic triangle with sides of length $s$ (opposite to the right angle)
and $r$, opposite to $\beta$.  By identity 2.6 we have $\beta=\sin^{-1}\big(\frac{\sinh(r)}{\sinh(s)}  \big)$.\\
{\bf 4.} Since the $M$-fibers $\{y\}\times M$ are orthogonal in $\cE_k(M)$ to the $\HH^k$-fibers $\HH^k\times\{ v\}$, items 1,2, and 3 above
imply that the $\bS^{k-1}$-fibers, the $\bS^{n-1}$-fibers and $(0,\pi/2)$-fibers are mutually orthogonal 
in $\bS^{k-1}\times\, \bS^{n-1}\times \big( 0\, ,\, \pi/2  \big)$ with the metric $\Xi_* f$.\\
{\bf 5.}
The map \begin{center}$\Xi'=(\Xi, s):E_k(M)\ra \bS^{k-1}\times\, \bS^{n-1}\times \big( 0\, ,\, \pi/2  \big)\times\R^+$\end{center} 
\noindent gives coordinates on $E_k(M)$.\\\\

\noindent {\bf \large  Section 4. Spherical Cuts on Hyperbolic Extensions.}

Let $(N^m, g)$ have center $o$.
Recall from the Introduction that the metric $g\0{r}$ on $\bS_r$ is called the {\it  spherical cut of $g$  at $r$}, and the metric $\hat{g}\0{r}\, =\, \big(\frac{1}{\sinh^2(r)}\big)\, g\0{r}$ is the {\it normalized spherical cut of $g$ at $r$.} 
\vspace{.0in}

Now let $(M^n,h)$ have center $o$. Thus we can write
$h=h_r+dr^2$, where each $h_r$ is a metric on $\bS^{n-1}$.
As before we denote by $f=\cE_k(h)$ the hyperbolic extension 
of $h$, and we write $f=f_s+ds^2$ on $\cE_k(M)-\{ o\}$;
each $f_s$ is a metric on $\bS^{n+k-1}$.
We use the map $\Xi=\Xi_s$ of Section 3 that gives coordinates on $S_s\big( \cE_k(M) \big)$.
Note that the metric  $\Xi_* f_s$ is a metric on $\bS^{k-1}\times\bS^{n-1}\times (0,\pi/2)$, and it 
is the expression of $f_s$ in the $\Xi$-coordinates.
\vspace{.1in}

\noindent {\bf Proposition 4.1.} {\it The expression of $f_s$ in the $\Xi$-coordinates
is given by} \begin{center}$\Xi_* f_s\,=\,\Big( \sinh^2(s)  \, \cos^2\, (\beta)\,\Big)\, \sigma\0{\bS^{k-1}}\, \, +\,\,
h_r\,\,+\,\,\Big( \,\sinh^2(s)  \,\Big)\,d\beta^2$\end{center}
\noindent {\it where $r=\sinh^{-1}(\sinh (s)\, \sin (\beta))$  (see identity 2.6).}\vspace{.1in}

\noindent {\bf Remark 4.2.}  Note that the function $r=r(s,\beta)$ is the same function used in the Introduction for
the $\theta$-reparametrizations $\lambda=\lambda(\lambda',\theta)$.
\vspace{.1in}

\noindent {\bf Proof.} By Remark 3.1(4) we have that $\Xi_* f_s$ has the form $A+B+C$, where $A(u,\beta)$ is a metric on
$\bS^{k-1}\times\{ u\}\times \{ \beta\}$, $B(w,\beta)$ is a metric on
$\{ w\}\times\bS^{n-1}\times \{ \beta\}$ and $C(u,\beta)$ is a metric on
$\{w\}\times \{ u\}\times (0,\pi/2)$, i.e. $C=f(w,u,\beta)\,d\beta^2$, for some positive function $f$.\vspace{.1in}

Now, by definition we have 
\begin{center}$f\,=\,\cosh ^2(r)\, \sigma\0{\HH^k}\, +\, h_r+dr^2\, =\,\cosh ^2(r)\bigg( \sinh^2(t)\sigma\0{\bS^{k-1}}+dt^2 \bigg)\, +\, h_r\, +\, dr^2 $\end{center}
\noindent By Proposition 2.7 and the identity $\cosh (r)\, \sinh(t)= \sinh(s)\, \cos\,(\beta)$
(which follows from the law of sines and the second law of cosines, also see identity 2.6) we can write
\begin{center}$f_s+ds^2=f\,=\,\big( \sinh^2(s)  \, \cos^2\, (\beta)\,\big) \sigma\0{\bS^{k-1}}\, \, +\,\,
h_r\,\,+\,\,\big( \sinh^2(s)  \big)\,d\beta^2\,\,+\,\, ds^2$\end{center}
\noindent This proves the proposition.\vspace{.1in}

Hence Proposition 4.1 gives the expression of
the spherical cut, at $s$, of the metric $f=\cE_k(h)$ in the $\Xi$-coordinates. The next corollary
does the same for the normalized spherical cut $\hat{f}$ of
$f$ at $s$.\vspace{.1in}

\noindent {\bf Corollary 4.3.} {\it The expression of ${\hat{f}}_s$ in the $\Xi$-coordinates
is given by} \begin{center}$\Xi_*\big(\,\, {\hat{f}}_s\,\,\big)\,=\, \cos^2\, (\beta)\, \sigma\0{\bS^{k-1}}\, \, +\,\,
\sin^2\, (\beta)\,\hh\0r\,\,+\,\,\,d\beta^2$\end{center}
\noindent {\it where $r$ as in Proposition 4.1.}
\vspace{.1in}

\noindent {\bf Proof.} Since $\sinh^2(r)\,\hh\0r=h_r$,
and $\sinh^2(s)\,{\hat{f}}_s=f_s$,
the corollary follows from Proposition 4.1 and identity 2.6.\\\\

\noindent {\bf \large  Section 5. Cut Limits and Proof of The Main Theorem.}

First a bit of notation.
Let $(N^m, g)$ have center $o$. Recall that we
can write the metric on $N-\{ o\}=\bS^{m-1}\times\R^+$ as \,\,$g=g\0{r}+dr^2$,
where $r$ is the distance to $o$. Let $A\sbs\bS^{m-1}$ be open and denote by $\rC A$ the open cone $A\times\R^+\sbs\bS^{m-1}\times\R^+\sbs M$. We write $A_r=\rC A\cap\bS_r(M)=A\times\{r\}$. We say that $\{ g_\ssl\}\0{\lambda}$ is an $\odot${\it-family of metrics over $A$} if each $g_\ssl$ is 
a metric defined on $\rC A$ and $g_\ssl$, and it can be written
in the form $g_\ssl=\big(g_\ssl\big)_r+dr^2$ on $\rC A$. 
We say that the $\{g_\ssl\}$ has {\it cut limit over $A$, at $b$,} if there is a 
$C^2$ metric $\hg_{_{\infty}}^b$ on $A$ such that (1.3) holds, where
the arrow in (1.3)  now means uniform convergence in the $C^2(A)$-norm on the 
space of $C^2$ metrics on  $A\sbs\bS^{m-1}$.
Also, {\it cut limits  over $A$, on $I$},
and {\it cut limits over $A$} are defined similarly.\\
%\vspace{.1in}

%\noindent{\bf Remark 5.1.} In the case of
%cut limits on an interval $I$ (including $I=\R$) we demand that the convergence
%in (1.3) is uniform in $b$ and uniform on $A$
%(see definition in paragraph before
%Remark 1.5).
%the $C^2$-convergence is uniform $C^2$-convergence
%with compact supports
%in the $I$-direction and uniform $C^2$-convergence over $A$ in the $\bS^{n-1}$-direction.
%\vspace{.1in}

 Let $M^n$ have metric $h$ and center $o$. As always we
identify $M-\{o\}$ with $\bS^{n-1}\times\R^+$ and 
$M$ with $\R^n$. Choose a center $o\in\HH^k\sbs\cE_k(M)$. Let $\{ h_\ssl\}_\ssl$ be a $\odot$-family of metrics on $M$, thus $o$ is a center for all $h_\ssl $. 
Denote by $f_\ssl=\cE_k(h_\ssl)$ the hyperbolic extension of
$h_\ssl$. We have that $\{f_\ssl\}_\ssl$ is a $\odot$-family on $\cE_k(M)$. %(see 3.3 \cite{O3}).
From now on we assume  $\theta \in (0,\pi/2]$ fixed. Next $\theta$-reparametrize $\{f_\ssl\}_\ssl$, that is, we use the change of
variables $\lambda=\lambda (\lambda')=\sinh^{-1}(\sinh(\lambda')\,\sin\,\theta)$. (Note that $\lambda'$ plays the role of the variable $s$
in identity 2.6, and $\lambda$ plays the role of $r$.)
We obtain in this way the $\odot$-family 
$\{f\0{\lambda(\lambda')}\}\0{\lambda'}$.
Write $S=\bS^{n+k-1}-\{\bS^{k-1}\coprod\bS^{n-1}\}$,
where $\bS^{k-1}\sbs\HH^k\times\{o\}$ and $\bS^{n-1}\sbs\{o\}\times M$.
%In the next proposition we fix $c\in \R$. Also we use concept of ``hyperbolic around the origin"
%given in the Introduction.
\vspace{.1in}

\noindent {\bf Proposition 5.1.}
{\it Assume that $\{h_\ssl\}$ has cut limits on the interval
$J_c=(-\infty, c]$, and that it is hyperbolic around the
origin. Then for each $c'<c+\ln\,\sin(\theta)$ the family  $\{f\0{\lambda(\lambda')}\}\0{\lambda'}$
has cut limits on $J_{c'}$ over $S$.}
\vspace{.1in}

\noindent {\bf Proof.} By hypothesis
$\{h_\ssl\}$  is hyperbolic around the
origin. Hence there is $B$ such that
\begin{equation*}
{\widehat{(h_\ssl)}}_{_{\lambda+b}}\,=\, \sigma\0{\bS^{n-1}}\,\,\,\,\,\,\,\,\,\,\,{\mbox{for all}}\,\,\,\,\,\,\,\,\,\,b\leq B
\tag{1}
\end{equation*}
Hence the metrics $h_\ssl$ are canonically hyperbolic 
on the ball of radius $\lambda+B$. 
Also, since we are assuming $\{h_\ssl\}$ has cut limits on
$J_c$ we have that
\begin{equation*}b\in J_c\,\,\,\,\,\,\,\,\,\,\,\,\,\,\,\,\,\,\Longrightarrow\,\,\,\,\,\,\,\,\,\,\,\,\,\,\,\,\,\,\,
{\widehat{(h_\ssl)}}_{_{\lambda+b}}\stackrel{C^2}{\longrightarrow}
\hh\0{\infty}^b\,\,\,\,\,\,\,{\mbox{as}}\,\,\,\,\lambda\ra\infty
\tag{2}
\end{equation*}
\noindent  uniformly on $\bS^{n-1}$ and uniformly with compact supports in the variable $b\in J_c$.\vspace{.1in}

As mentioned before we can write $f_\ssl\,=\,\big(f_\ssl\big)_s +ds^2$. We have to compute the limit of
${\widehat{(f\0{\lambda(\lambda')})}}\0{\lambda'+b}$,
as $\lambda'\ra\infty$.\vspace{.1in}
Let the $\Xi$-coordinates be as defined in Section 3 
for the space $(\cE_k(M),f)$.  From Corollary 4.3. we can express $\big(\hf_\ssl\big)_s$ in $\Xi$-coordinates:
\begin{center}$\Xi_*\big(\,\, {\widehat{(f\0{\lambda(\lambda')})}}\0{\lambda'+b}\,\,\big)\,=\, \cos^2\, (\beta)\, \sigma\0{\bS^{k-1}}\, \, +\,\,
\sin^2\, (\beta)\,{\widehat{(h\0{\lambda(\lambda')})}}\0{r(\lambda'+b,\beta)}\,\,+\,\,\,d\beta^2$\end{center}
\noindent  where $r=r(s,\beta)$ is given by identity 2.6 (see also  Proposition 4.1 and Remark 4.2). Therefore 
we want to find the limit of ${\widehat{(h\0{\lambda(\lambda')})}}\0{r(\lambda'+b,\beta)}$ as $\lambda'\ra\infty$. 
To do this take the inverse
of $\lambda=\lambda(\lambda')$, and we get $\lambda'=\lambda'(\lambda)=\sinh^{-1}\big( \frac{\sinh(\lambda)}{ \sin (\theta)}  \big)$.
Hence
\begin{equation}\displaystyle \lim_{\lambda'\ra\infty}{\widehat{(h\0{\lambda(\lambda')})}}\0{r(\lambda'+b,\beta)}\,\,=\,\, 
\lim_{\lambda\ra\infty}{\widehat{(h\0{\lambda})}}\0{\vartheta(\lambda,\beta, b)},\tag{3}\end{equation}
\noindent where 
\begin{center}$\vartheta(\lambda,\beta, b)=r\big(\lambda'(\lambda)+b,\beta\big)=
\sinh^{-1}\bigg(\,\sinh\,\bigg\{ \,\,b+ \sinh^{-1}\bigg( \frac{\sinh(\lambda)}{\sin(\theta)} \bigg) \,\,\,\bigg\}\,\,\sin\, (\beta)\,    \bigg)$\end{center}
\noindent and a straightforward calculation shows 
\begin{equation*}
\lim_{\lambda\ra\infty}\Big(\,\vartheta (\lambda,\beta, b)\,-\, \lambda\,\Big)\,=\, b\,+\, \ln\big( \,\frac{\sin(\beta)}{\sin(\theta)}\,\big).
\tag{4}
\end{equation*}
\noindent This convergence is 
uniform with compact supports in the $C^2(S)$-topology (see caveat below).
Choose $c'\in \R$ such that $c'<c-\ln\big(\frac{\sin(\pi/2) }{\sin(\theta)}\big)=c+\ln\,\sin (\theta)$.
Since $\beta\in(0,\pi/2)$  we get 
\begin{equation}
b\,\,\in\,\, J_{c'}\,\,\,\,\,\,\,\,\,\,\Longrightarrow\,\,\,\,\,\,\,\,\,\,\big(\, b\,+\, \ln\big( \,\frac{\sin(\beta)}{\sin(\theta)}\,\big)\,\big)\,\,\in \,\,J_c.
\tag{5}\end{equation}
Hence from (2), (3), (4)  and (5) we get 
\begin{equation}\displaystyle\lim_{\lambda'\ra\infty}{\widehat{(h\0{\lambda(\lambda')})}}\0{r(\lambda'+b,\beta)}\,\,=\,\,\hh\0{\infty}
^{b+\ln\,(\frac{\sin\,\beta}{\sin\,\theta})}\tag{6}\end{equation}
\noindent {\it Caveat. The limit (3) (hence also in (6)) is uniform} {\sf with compact supports} {\it in the $\beta$ direction, but not uniform in the
$\beta$ direction. The problem occurs when $\beta \ra 0$.}
\vspace{.1in}

We next deal with the problem mentioned
in the caveat; that is, we have to show that the limit in (6) in uniform in the variable
$\beta\in (0,\pi/2)$ (not just uniform with compact supports). The
convergence in (4) (hence in (6)) is uniform for $\beta$ near $\pi/2$, but the convergence in (4) is certainly
not uniform near 0. Here is where we will need the extra condition of the family being
hyperbolic near the origin.  We will need the following claim.\vspace{.1in}

\noindent{\bf Claim.} {\it  Let $c, B, \theta\in \R$. Choose $c'$ with
$c'<c+\ln\,\sin\,\theta$.
Then there is $\beta_1>0$
such that $r(\lambda'+c',\beta_1)\leq
\lambda(\lambda')+ B$, for every $\lambda'$ sufficiently large.}\vspace{.1in}

\noindent{\bf Proof of the claim.} A calculation
shows that taking
$\beta_1=\sin^{-1}(e^{2(B-c-1})$ works.
(Find the limit $\lambda'\ra\infty$ of both
terms in the inequality, and use the fact that
$c'<c+\ln\,\sin\,\theta$.) This proves the claim.
\vspace{.1in}

Since the function $r=r(s,\beta)$ is increasing in
both variables, the claim implies that
$r(\lambda'+b,\beta)\leq
\lambda(\lambda')+ B$, for every $b\leq c'$,
$\beta\leq \beta_1$ and $\lambda'$ sufficiently
large (how large not depending on $b$, nor $\beta$).
This together with (1) imply that
for every $b\leq c'$,
$\beta\leq \beta_1$ and $\lambda'$ sufficiently
large we have
\begin{center}${\widehat{(h\0{\lambda(\lambda')})}}\0{r(\lambda'+b,\beta)}\,\,=\,\,\sigma\0{\bS^{n-1}}$\end{center}
\noindent Hence for every $b\in J_{c'}$ and
$\beta\leq \beta_1$ we have
\begin{center}$\displaystyle\lim_{\lambda'\ra\infty}{\widehat{(h\0{\lambda(\lambda')})}}\0{r(\lambda'+b,\beta)}\,\,=\,\,\sigma\0{\bS^{n-1}}$\end{center}
Since $\beta_1>0$ the problem mentioned in the caveat (i.e. when $\beta\ra0$) has been removed.  This proves the proposition.\vspace{.1in}

Taking $c\ra\infty$ in Proposition 5.1 gives the
following corollary.\vspace{.1in}

\noindent {\bf Corollary 5.2.}
{\it Assume that $\{h_\ssl\}$ has cut limits, and that it is hyperbolic around the
origin. Then $\{f\0{\lambda(\lambda')}\}\0{\lambda'}$
has cut limits over $S$.}
\vspace{.1in}

\noindent {\bf Proof of the Main Theorem.}
Note that the only difference between Corollary
5.2 and the Main Theorem is that in the
corollary the cut limits exist {\sf over} $S\sbs\bS^{n+k-1}$. Hence
we have to show that the existence of cut limits
over $S$ implies the existence of cut limits on the
whole of $\bS^{n+k-1}$. Corollary 5.2 and 1.3
in the Introduction imply\vspace{.1in}

\noindent \hspace{2in}$
\big|\,\,{\widehat{\big((f_\ssl)|_S\big)}}_{_{\lambda'+b}} \,\,-\,\, \hf_{_{\infty}}^b\,\,\big|_{C^2(S)}\,\,\longrightarrow\,\, 0
$  \,\,\,\,\,\,\,as\,\,\,\,\,\, $\lambda'\ra \infty$\vspace{.1in}

\noindent where $\hf_{_{\infty}}^b$ is a metric on $S$. In particular for every $b$ the one-parameter family
 ${\widehat{\big((f_\ssl)|_S\big)}}_{_{\lambda'+b}}$
is Cauchy, that is \vspace{.1in}
\begin{equation}
\big|\,\,{\widehat{\big((f\0{\lambda(\lambda'_1)})|_S\big)}}_{_{\lambda'_1+b}} \,\,-\,\, {\widehat{\big((f\0{\lambda(\lambda'_2)})|_S\big)}}_{_{\lambda'_2+b}}\,\,\big|_{C^2(S)}\,\,\longrightarrow\,\, 0  \tag{7}\end{equation}
\noindent uniformly on $S$ as $\lambda'_1$, $\lambda'_2$ $\ra \infty$.  But since $S$ is dense in $\bS^{n+k-1}$ 
we get that  $|\,g|\0{S}\,|_{C^2(S)}=|g|_{C^2(\bS^{n+k-1})}$, for any $C^2$ (pointwise) bilinear form $g$ on $\bS^{n+k-1}$.
Therefore
we can drop the restriction ``$|_S$" in 
(7) to get
\begin{center}
$\big|\,\,{\widehat{\big(f\0{\lambda(\lambda'_1)}\big)}}_{_{\lambda'_1+b}} \,\,-\,\, {\widehat{\big(f\0{\lambda(\lambda'_2)}\big)}}_{_{\lambda'_2+b}}\,\,\big|_{C^2(\bS^{n+k-1})}\,\,\longrightarrow\,\, 0 
\,\,\,\,\,{\mbox{as}}\,\,\,\,\,\lambda'\ra 0$ \end{center}
\noindent This implies that the family
${\widehat{\big(f_\ssl\big)}}_{_{\lambda'+b}}$
is Cauchy. Since the space of $C^2$ metrics on $\bS^{n+k-1}$ with the $C^2$ norm is a complete metric space
the Cauchy sequence above converges to some $\hf\0{\infty}^b$. Note that $\hf\0{\infty}^b$
is a symmetric bilinear form on $\bS^{n+k-1}$,
and it is positive definite on $S$. It remains to prove
that $\hf_{_{\infty}}^b$ is also positive definite
outside $S$. Recall $S=\bS^{n+k-1}-(\bS^{k-1}\coprod\bS^{n-1})$.
But it is straightforward to verify that  
we have $\hf\0{\infty}^b|_{\bS^{k-1}}=\sigma\0{\bS^{k-1}}+\sigma\0{\HH^n}$. On the other hand on $\bS^{n-1}$ we have $\beta=\pi/2$,
hence $\lambda=\lambda'$. Also by
definition we have $f\0{\lambda}=
\cosh ^2(r)\sigma\0{\HH^k}+h_\lambda$. But on $M_o$ we get $r=s$. Therefore
$${\widehat{\big((f_\ssl)|_{\bS^{n-1}}\big)}}\0{\lambda'+b}=
{\widehat{\big((f_\ssl)|_{\bS^{n-1}}\big)}}_{_{\lambda+b}}=
{\mbox{cotanh}}^2(\lambda +b)\sigma\0{\HH^k}+
{\widehat{(h\0{\lambda})}}\0{\lambda+b}
\longrightarrow
{\mbox{{\mbox{cotanh}}}}^2(\lambda +b)\sigma\0{\HH^k}+\hh\0{\infty}^b$$
\noindent Consequently $\hf\0{b+\infty}$
is positive definite on $\bS^{n-1}$. Thus it is
positive definite outside $S$. This proves the 
Main Theorem.

Pedro Ontaneda

SUNY, Binghamton, N.Y., 13902, U.S.A.

\end{document}